\documentclass{sigcomm-alternate}
\usepackage{amsfonts,latexsym,amsmath,amstext,amssymb,verbatim,epsfig,psfrag}
\usepackage{colordvi,pstcol}
\usepackage{graphicx}
\usepackage{psfrag}

\newcommand\II{\mathcal{I}}

\def\R{\mbox{I\hspace{-.15em}R}}

\def\Pb{\bold{P}}

\def\Hb{\bold{H}}

\def\Ib{\bold{I}}
\def\Jb{\bold{J}}

\begin{document}

\title{Introducing One Step Back Iterative Approach to Solve Linear and Non Linear Fixed Point Problem}

\numberofauthors{1}
\author{
   \alignauthor{Dohy Hong}\\
   \affaddr{Alcatel-Lucent Bell Labs}\\
   \affaddr{Route de Villejust}\\
   \affaddr{91620 Nozay, France}\\
   \email{\normalsize dohy.hong@alcatel-lucent.com}
}

\maketitle

\begin{abstract}
In this paper, we introduce a new iterative method which we call one step back approach: the main idea is
to anticipate the consequence of the iterative computation per coordinate and to optimize on the choice
of the sequence of the coordinates on which the iterative update computations are done.
The method requires the increase of the size of the state vectors and one iteration step loss
from the initial vector.
We illustrate the approach in linear and non linear iterative equations.
\end{abstract}


\begin{psfrags}
\section{Introduction}\label{sec:intro}
Iterative methods to solve large sparse systems have been gaining interests in very different research areas
and a large number of approaches have been studied: starting from Jacobi or Gauss-Seidel iteration
methods, more recent works focuses on relaxation methods, Krylov subspace methods, the use of
preconditioners, matrix decomposition, parallel computation etc. to solve linear or non linear
problems
\cite{Saad}, \cite{cg}, \cite{arnoldi}, \cite{greenbaum}, \cite{stewart},
\cite{gantmacher2000theory}, \cite{berman1994nonnegative}, \cite{varga2009} to cite
few of them.

In this paper, we propose a new iterative algorithm based on the anticipated consequence
of the iteration at each coordinate level: 
we believe that this approach may bring significant improvement in a large class of linear and
non linear problems and it seems that such an approach have not yet been studied by the
research community.
More precisely, we study the computation of the fixed point $X\in \R^N$ solving 
$$ 
X = \Hb (X)
$$
starting from the initial condition $X_0 \in \R^N$.

In linear algebra, this approach is equivalent to the idea
of fluid diffusion that was called D-iteration in \cite{d-algo}.

In this paper, we will not consider the convergence issue and assume that the Jacobi or the
asynchronous version of Gauss-Seidel style iteration of $X_{n+1} = \Hb(X_n)$ 
converges to a unique fixed point (that's basically contracting operators, cf. for instance \cite{Baudet, bertsekas, frommer}).

In Section \ref{sec:algo}, we define the notations and the proposed One Step Back (OSB)
algorithm.
Section \ref{sec:illustration} illustrates two concrete examples to solve linear
and non linear systems. 
\section{OSB algorithm}\label{sec:algo}
We will use the following notations:
\begin{itemize}
\item $\Ib \in \R^{N\times N}$ the identity matrix;
\item $\Jb_i$ the matrix with all entries equal to zero except for
  the $i$-th diagonal term: $(\Jb_i)_{ii} = 1$;
\item $\Hb : \R^{N} \to \R^{N}$ an operator;
\item $\II =\{i_1,...,i_n,...\}$ a sequence of the coordinates (or nodes), $i_n \in \{1,..,N\}$.
\end{itemize}

\subsection{Iterative equations}
The OSB iteration is defined by the triplet $(\Hb,X_0,\II)$ and exploits two state vectors of size $N$,
$H_n$ (history) and $F_n$ (residual fluid) based on the following iterative equations:
\begin{eqnarray}
H_0 &=& X_0\nonumber\\
H_n &=& H_{n-1} + \Jb_{i_n} (F_{n-1})\label{eq:defH}
\end{eqnarray}
and
\begin{eqnarray}
F_0 &=& \Hb(H_0) - H_0\nonumber\\
F_n &=& (\Ib - \Jb_{i_n})(F_{n-1}) + \Hb(H_n) - \Hb(H_{n-1}).\label{eq:defF}
\end{eqnarray}

\subsection{Sequence choice and distributed computation}
The choice of the sequence $\II$ depends on the computation costs structure that need
to be optimized. One may consider those introduced in \cite{dohy}, \cite{revisit} or \cite{diff}
depending on the context.

For the distributed computation architecture, one may also use the one proposed in
\cite{part} for the linear case: the architecture for linear or non linear cases
should be a priori the same and benefits from the asynchronous properties of
the proposed method.

\section{Illustration}\label{sec:illustration}

\subsection{Linear equation: D-iteration}
If we take $\Hb = \Pb$ a linear operator, we have:
\begin{eqnarray*}
F_n &=& (\Ib - \Jb_{i_n})(F_{n-1}) + \Hb(H_n) - \Hb(H_{n-1})\\
 &=& (\Ib - \Jb_{i_n})F_{n-1} + \Pb \Jb_{i_n} F_{n-1}.
\end{eqnarray*}
Then we find back D-iteration equation cf. \cite{d-algo}.
The computation gain for this case has been studied in details for instance in
\cite{dohy, part, revisit} and is skipped here.

\subsection{Non linear equation: simple example}
Consider the non-linear fixed point problem of dimension three:
$$
\Hb(x,y,z) = (\sqrt{xy}+1, (x+z)/4+1, (x+y)/4).
$$
The fixed point of the above equation is $((23+\sqrt{379})/10, (23+\sqrt{379})/30+16/15, (23+\sqrt{379})/30+4/15) \sim (4.247, 2.482, 1.682)$.
The updates of the fluid vector $F_n$ is based on the three increment functions:
\begin{eqnarray*}
\Hb(x',y,z)-\Hb(x,y,z) &=& (\sqrt{y}(\sqrt{x'}-\sqrt{x}), \frac{(x'-x)}{4}, \frac{(x'-x)}{4})\\
\Hb(x,y',z)-\Hb(x,y,z) &=& (\sqrt{x}(\sqrt{y'}-\sqrt{y}), 0, \frac{(y'-y)}{4})\\
\Hb(x,y,z')-\Hb(x,y,z) &=& (0, \frac{(z'-z)}{4}, 0).
\end{eqnarray*}

Figure \ref{fig:exm} shows the convergence of Jacobi, Gauss-Seidel and OSB iterations. One iteration is here
defined as an update of three entries (with multiplicity for OSB) of the iterated vector.
For OSB, we choose the $n$-th coordinate equal to the $\arg\max$ of $F_n$ in absolute  value.
The initial condition is set to $(4.2, 1, 1.5)$.

\begin{figure}[htbp]
\centering
\includegraphics[angle=-90, width=7cm]{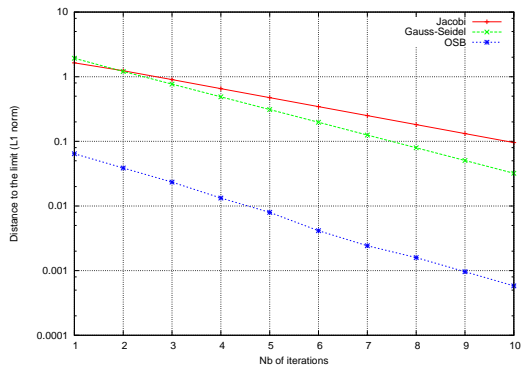}
\caption{Illustration: $N=3$.}
\label{fig:exm}
\end{figure}

Figure \ref{fig:exm} illustrates
well the advantage of optimizing the coordinate sequence order when the initial vector has some coordinates
more closer to the limit than others: the gain factor compared to Jacobi iteration is here of two orders of magnitude
at iteration 10!
Note that in this case, starting from $(0,0,0)$, there is no big differences between those three methods
(the matrix size is too small to observe significant differences). 

This example is only shown for the sole purpose of the illustration of the potential impact
and no theoretical guarantee on the convergence gain is given here. However, the author believes that
we should have cases where the OSB method may provide substantial convergence improvement, since it
can be applied to a very wide range of fixed point problems. 
The first factor to be considered to determine whether OSB method can improve the iterative computation cost
is the complexity of the computation of the increment $\Hb(H_n) - \Hb(H_{n-1})$, which may simplify or
introduce a cost overhead to Gauss-Seidel style {\em normal} iteration, and to compare this complexity
to the gain brought by coordinate level optimization.

Note finally that the one step iteration computation lost at the first iteration is recovered at the end
considering $H_n+F_n$ as the estimator (instead of $H_n$).

\section{Conclusion}\label{sec:conclusion}
In this paper, we described a new iterative method and illustrated its applications
to two simple fixed point problems. 

\end{psfrags}
\bibliographystyle{abbrv}
\bibliography{sigproc}

\end{document}